\pdfoutput=1
\RequirePackage{ifpdf}
\ifpdf % We~are running pdfTeX in pdf mode
\documentclass[pdftex]{sigma}
\else
\documentclass{sigma}
\fi

\numberwithin{equation}{section}
\newtheorem{Theorem}{Theorem}[section]
{\theoremstyle{definition}
\newtheorem{Example}[Theorem]{Example}}

\newcommand*{\mM}{\mathbb M} \newcommand*{\mR}{\mathbb R}
\newcommand*{\mS}{\mathbb S} \newcommand*{\mT}{\mathbb T}
\newcommand*{\mZ}{\mathbb Z}

\newcommand*{\fI}{\mathfrak I} \newcommand*{\fJ}{\mathfrak J}
\newcommand*{\fK}{\mathfrak K} \newcommand*{\fL}{\mathfrak L}

\newcommand*{\rmd}{\mathrm d} \newcommand*{\rme}{\mathrm e}
\newcommand*{\rmi}{\mathrm i}

\numberwithin{figure}{section}

\begin{document}
\allowdisplaybreaks

\newcommand{\arXivNumber}{2209.02620}

\renewcommand{\PaperNumber}{084}

\FirstPageHeading

\ShortArticleName{Three Examples of Dynamical Systems}

\ArticleName{Three Examples in the Dynamical Systems Theory}

\Author{Mikhail B.~SEVRYUK}

\AuthorNameForHeading{M.B.~Sevryuk}

\Address{V.L.~Talrose Institute for Energy Problems of Chemical Physics,\\
 N.N.~Sem\"enov Federal Research Center of Chemical Physics,\\ Russian Academy of Sciences, Leninski\u{\i} Prospect 38, Bld.~2, Moscow 119334, Russia}
\Email{\href{mailto:2421584@mail.ru}{2421584@mail.ru}}

\ArticleDates{Received September 10, 2022, in final form October 18, 2022; Published online October 29, 2022}

\Abstract{We present three explicit curious simple examples in the theory of dynamical systems. The first one is an example of two analytic diffeomorphisms $R$, $S$ of a closed two-dimensional annulus that possess the intersection property but their composition $RS$ does not ($R$ being just the rotation by $\pi/2$). The second example is that of a~non-Lagrangian $n$-torus $L_0$ in the cotangent bundle $T^\ast\mT^n$ of $\mT^n$ ($n\geq 2$) such that $L_0$ intersects neither its images under almost all the rotations of $T^\ast\mT^n$ nor the zero section of $T^\ast\mT^n$. The third example is that of two one-parameter families of analytic reversible autonomous ordinary differential equations of the form $\dot{x}=f(x,y)$, $\dot{y}=\mu g(x,y)$ in the closed upper half-plane $\{y\geq 0\}$ such that for each family, the corresponding phase portraits for $0<\mu<1$ and for $\mu>1$ are topologically non-equivalent. The first two examples are expounded within the general context of symplectic topology.}


\Classification{57R17; 53D12; 37C15}

\section{Introduction}\label{introduction}

It is common knowledge what a significant role examples and counterexamples (the difference between these two concepts is mainly psychological and historical, not objective) play in mathematics and especially in mathematical education. For instance, V.I.~Arnold, one of the greatest mathematicians of the twentieth century, wrote that ``no matter how much time one could save using the deductive methods (`from general to specific'), the value of a lecture for a student consists of merely a number of well-explained and thoroughly understood examples'' \cite{Arnold14}. He constantly emphasized the perniciousness of ``the axiomatic-deductive method'' which has ``led to the banishment of all examples (and especially the motives for introducing definitions) in the teaching of mathematics at every level''~\cite{Arnold96}.

I agree completely with S.Yu.~Yakovenko that while ``mathematical anecdotes mention great mathematicians whom examples only distracted from developing general theories'', ``Arnold was the opposite: examples were the alpha and omega of his approach'' \cite{Yakovenko14}.\footnote{I witnessed the following statement made by Arnold in his characteristic provocative style at his Moscow seminar on singularity theory: ``A million examples are more convincing than a proof''. Indeed, a series of carefully chosen examples can be more convincing than a lengthy and tedious formal proof. Of course, Arnold did not assert that examples can serve as a \emph{substitution} for a proof.}

Another outstanding mathematician of the last century, I.R.~Shafarevich, wrote that ``the meaning of a mathematical notion is by no means confined to its formal definition; in fact, it may be rather better expressed by a (generally fairly small) sample of the basic examples, which serve the mathematician as the motivation and the substantive definition, and at the same time as the real meaning of the notion'' \cite{Shafarevich05}.

Even using only free Internet resources, it is easy to estimate the number of mathematical works whose titles contain the words ``example'' or ``counterexample''. For instance, as of September~1, 2022, the site \url{https://mathscinet.ams.org/mrlookup} gives us 9182 items with the word ``example/examples'' in the title and 2332 items with the word ``counterexample/counterexamples'' in the title. Besides, also as of September~1, 2022, the numbers of \emph{books} recorded in \url{https://zbmath.org/} whose titles contain the words ``example'', ``examples'', ``counterexample'', ``counterexamples'' are equal to 115, 656, 3, 33, respectively. Of course, some of these books are different editions (or even editions in different languages) of the same book, but, on the other hand, not all books are indexed in \texttt{zbMATH}; for instance, V.~Boss'\footnote{A pseudonym of V.I.~Opo\u{\i}tsev.} manual~\cite{Boss20}, widely known in Russia, is not displayed in \texttt{zbMATH}. This manual belongs to the sort of books (very popular in the mathematical community) that are collections of counterexamples in one or more branches of mathematics. I cannot resist mentioning the first (according to \texttt{zbMATH}) and probably the most famous book of this kind, namely, the masterpiece~\cite{GelbaumOlmsted03} whose first edition was published in 1964. As declared in the preface to \cite{GelbaumOlmsted03}, ``At the risk of oversimplification, we might say that (aside from definitions, statements, and hard work) mathematics consists of two classes~-- proofs and counterexamples''. The similar statement, ``Mathematics consists of two things~-- theorems and counterexamples'', is sometimes cited (for instance, this is the epigraph to the preface to the book~\cite{Boss20}) as being due to G.~Polya but without any precise reference. I~have also failed to find the source.

This note describes three rather simple but curious, from the author's viewpoint, examples in the theory of dynamical systems. The first two examples deal with symplectomorphisms and their generalizations whereas the third one concerns reversible vector fields.

\section{The intersection property in dimension two}\label{two}

Given an arbitrary $n$-dimensional manifold $L^n$, let $\lambda^1=p\,\rmd q$ be the standard action form (also known as the Liouville $1$-form, the Poincar\'e $1$-form, the canonical $1$-form, the tautological $1$-form, the distinguished $1$-form, or the symplectic potential) in the cotangent bundle $T^\ast L$, then $\rmd\lambda^1=\rmd p\wedge\rmd q$ is the canonical symplectic structure in $T^\ast L$. In the sequel, while considering $T^\ast L$ as a symplectic manifold, we will always assume that the symplectic structure in $T^\ast L$ is canonical. A symplectomorphism $A\colon T^\ast L\to T^\ast L$ is said to be \emph{exact} if the $1$-form $A^\ast\lambda^1-\lambda^1$ is not only closed but even exact. For instance, the phase flow mapping $\phi_0^1$ of any Hamiltonian (not necessarily autonomous) vector field $V_t$ in $T^\ast L$ is an exact symplectomorphism of $T^\ast L$. It is clear that the inverse map of an exact symplectomorphism $A$ of $T^\ast L$ is also an exact symplectomorphism and that the composition of two exact symplectomorphisms $A$, $B$ of $T^\ast L$ is again an exact symplectomorphism: if $A^\ast\lambda^1-\lambda^1=\rmd X$ and $B^\ast\lambda^1-\lambda^1=\rmd Y$ ($X,Y\colon T^\ast L\to\mR$) then
\begin{gather*}
\big(A^{-1}\big)^\ast\lambda^1-\lambda^1 = \big(A^{-1}\big)^\ast\big(\lambda^1-A^\ast\lambda^1\big) = -\big(A^{-1}\big)^\ast\rmd X = -\rmd\big(\big(A^{-1}\big)^\ast X\big), \\
(AB)^\ast\lambda^1-\lambda^1 = B^\ast A^\ast\lambda^1-B^\ast\lambda^1+B^\ast\lambda^1-\lambda^1 = B^\ast\rmd X+\rmd Y = \rmd(B^\ast X+Y).
\end{gather*}

Suppose that $n=1$ and $L=\mS^1=\mR/2\pi\mZ$ is the standard circle, so that $T^\ast L\cong\mR\times\mS^1$ is a cylinder (or an open annulus). One says that a diffeomorphism $A$ of $\mM^2=\mR\times\mS^1$ possesses the \emph{intersection property} (or \emph{self-intersection property}) if for any simple\footnote{That is, without self-intersections.} closed curve $\gamma\subset\mM^2$ homotopic to $\{0\}\times\mS^1$ \big($0\in\mR$, $\gamma\cong\mS^1$\big), the image of $\gamma$ under the action of $A$ intersects~$\gamma$.\footnote{One sometimes requires that $A\gamma\cap\gamma\neq\varnothing$ for any non-null-homotopic (i.e., not contractible) simple closed curve $\gamma$.} For instance, any exact symplectomorphism $A$ of $T^\ast\mS^1$ (e.g., any twisted rotation $(p,q)\mapsto(p,q+v(p))$ which is the phase flow mapping $\phi^1$ of the Hamiltonian vector field $v(p)\partial_q$ with the Hamilton function $\int_0^p v(\tau)\,\rmd\tau$) possesses the intersection property. Indeed, if a simple closed curve $\gamma\subset T^\ast\mS^1$ homotopic to $\{0\}\times\mS^1$ and the curve $A\gamma$ do not intersect then the (algebraic) area $\int_S\rmd\lambda^1$ of the domain $S$ bounded by $\gamma$ and $A\gamma$ is non-zero. On the other hand, this area is equal to
\[
\int_{\partial S}\lambda^1 = \int_{A\gamma}\lambda^1-\int_\gamma\lambda^1 = \int_\gamma\big(A^\ast\lambda^1-\lambda^1\big) = 0.
\]
Many dynamical features of exact symplectomorphisms of $T^\ast\mS^1$ (they are also called \emph{globally} area preserving mappings) are exhibited by arbitrary diffeomorphisms with the intersection property. In particular, for such diffeomorphisms, one can develop the KAM (Kolmogorov--Arnold--Moser) theory, see, e.g., the works \cite{Moser62,Moser01,Russmann83,SiegelMoser95}; for some recent results and generalizations see, e.g., the papers \cite{HuangLiLiu18,HuangLiLiu22,MaXu22,XiaYu22}.

It is obvious that if a diffeomorphism of $\mM^2$ possesses the intersection property, so does the inverse diffeomorphism. However, the composition of two diffeomorphisms of $\mM^2$ possessing the intersection property does not necessarily possess the intersection property: the diffeomorphisms of $\mM^2$ with the intersection property do not constitute a group.

Let $(x,y)$ be Cartesian coordinates in $\mR^2$.

\begin{Example}\label{exampletwo}
Two analytic diffeomorphisms $R,S\colon M\to M$ of the closed annulus $M = \bigl\{(x,y) \bigm| 1\leq x^2+y^2\leq 25\bigr\}$ with the intersection property such that $R$ is the rotation $R\colon (x,y)\mapsto(-y,x)$ by $\pi/2$, $S$ is arbitrarily close to $R$, and the composition $RS$ does \textit{not} possess the intersection property.

In the annulus $M$, introduce the new coordinate frame $(\rho, \vartheta\bmod 2\pi)$ with $1\leq\rho\leq 5$ according to the formulas
\[
x=\chi(\rho)\cos\vartheta, \qquad y=\eta(\rho)\sin\vartheta, \qquad
\chi(\rho)=\frac{7\rho-5}{\rho+1}, \qquad \eta(\rho)=\frac{\rho+5}{7-\rho}.
\]
This coordinate frame is well defined since $\chi(1)=\eta(1)=1$, $\chi(5)=\eta(5)=5$, and $\chi'(\rho)=12/(\rho+1)^2>0$, $\eta'(\rho)=12/(7-\rho)^2>0$ everywhere. One has $\big\{x^2+y^2=1\big\} = \{\rho=1\}$, $\big\{x^2+y^2=25\big\} = \{\rho=5\}$, and it is very easy to see that $1<\eta(\rho)<\rho<\chi(\rho)<5$ for $1<\rho<5$. In particular, $\chi(3)=4$ and $\eta(3)=2$. In the sequel, we will also use the equalities $\chi(7/5)=2$ and $\eta(23/5)=4$.

Let $P(\rho) = (\rho-1)(\rho-3)^2(\rho-5)$. Then $P(1)=P(3)=P(5)=0$ and $P(\rho)<0$ for $1<\rho<3$ and $3<\rho<5$. Also it is very easy to verify that the global minimum of the derivative $P'(\rho) = 4(\rho-3)\big(\rho^2-6\rho+7\big)$ over the interval $1\leq\rho\leq 5$ is equal to $P'(1)=-16$. In fact, since $P''(\rho) = 4\big(3\rho^2-18\rho+25\big)$, the polynomial $P'(\rho)$ also attains a local minimum at $\rho=3+(2/3)^{1/2}$ in this interval, but this minimum is equal to $-(16/3)(2/3)^{1/2} = -4.35\ldots$ Thus, for $0<\varepsilon_1<1/16$, the mapping $\rho\mapsto\rho+\varepsilon_1P(\rho)$ is a diffeomorphism of the interval $1\leq\rho\leq 5$ that leaves the points $\rho=1,3,5$ fixed and moves all the other points towards the left.

Let also $\Theta(\vartheta) = \sin 2\vartheta$. Since the derivative $\Theta'(\vartheta) = 2\cos 2\vartheta$ nowhere exceeds $2$ in absolute value, the mapping $\vartheta\mapsto\vartheta+\varepsilon_2\Theta(\vartheta)$ with $0<\varepsilon_2<1/2$ is a diffeomorphism of the circle $\mS^1$ with fixed points $0$, $\pi$ (repellers) and $\pi/2$, $3\pi/2$ (attractors), all the other points of the circle being shifted towards these attractors.

Now consider two ellipses $\gamma = \big\{x^2/4+y^2/16=1\big\}$ and $R\gamma = \big\{x^2/16+y^2/4=1\big\} = \{\rho=3\}$ in~$M$ and the analytic diffeomorphism
\[
G\colon \ (\rho, \vartheta) \mapsto \bigl(\rho+\varepsilon_1P(\rho), \vartheta+\varepsilon_2\Theta(\vartheta)\bigr)
\]
of the annulus $M$ ($0<\varepsilon_1<1/16$, $0<\varepsilon_2<1/2$). This diffeomorphism leaves the ellipse $R\gamma$ invariant (since $P(3)=0$) and takes the ellipse $\gamma$ to a closed curve $\delta$ lying strictly inside $\gamma$. Indeed, the four vertices of $\gamma$ are $(\rho=7/5, \vartheta=0)$, $(\rho=7/5, \vartheta=\pi)$, $(\rho=23/5, \vartheta=\pi/2)$, and $(\rho=23/5, \vartheta=3\pi/2)$. The diffeomorphism $G_1\colon (\rho, \vartheta) \mapsto \bigl(\rho, \vartheta+\varepsilon_2\Theta(\vartheta)\bigr)$ leaves these vertices invariant and moves \big(``along'' the ellipses $\big\{\rho=\rho^0\big\}$, $7/5<\rho^0<23/5$\big) all the other points of $\gamma$ to some points strictly inside $\gamma$. The diffeomorphism $G_2\colon (\rho, \vartheta) \mapsto \bigl(\rho+\varepsilon_1P(\rho), \vartheta\bigr)$ moves the vertices of $\gamma$ to points strictly inside $\gamma$ because $P(7/5)<0$ and $P(23/5)<0$. On the other hand, if a point $C=(x,y)$ lies strictly inside $\gamma$, so does the point $(\tilde{x},\tilde{y})=G_2C$ because $|\tilde{x}|\leq |x|$ and $|\tilde{y}|\leq |y|$. Consequently, $\delta=G\gamma=G_2G_1\gamma$ lies strictly inside $\gamma$ (see Figure~\ref{fig21}).

\begin{figure}[t]\centering
\includegraphics{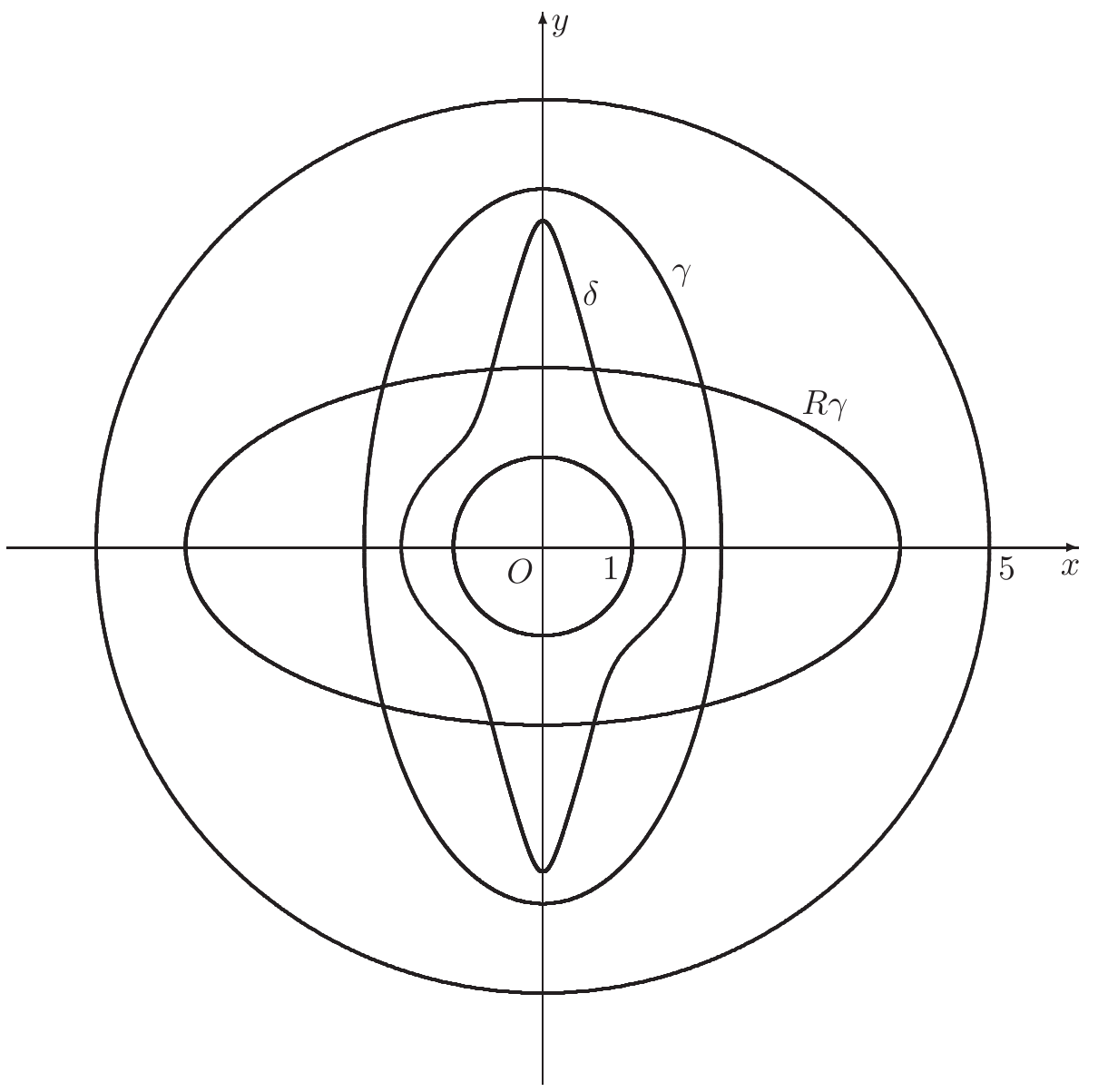}
\caption{The key curves in Example~\ref{exampletwo} for $\varepsilon_1=0.05$ and $\varepsilon_2=0.4$.}\label{fig21}
\end{figure}

The diffeomorphism $S=GRG^{-1}\colon M\to M$ possesses the intersection property (since $R$ does) and takes $R\gamma$ to $\delta$ and vice versa:
\begin{gather*}
SR\gamma = GRG^{-1}R\gamma = GRR\gamma = G\gamma = \delta, \\
S\delta = GRG^{-1}G\gamma = GR\gamma = R\gamma.
\end{gather*}
The composition $RS$ takes the curve $R\gamma$ to the curve $R\delta$ which lies strictly inside $R\gamma$ and takes the curve $\delta$ to the curve $RR\gamma=\gamma$ which lies strictly outside $\delta$. Thus, $RS$ does not possess the intersection property.

By choosing $\varepsilon_1$ and $\varepsilon_2$ to be arbitrarily small, one can make the diffeomorphism $G$ arbitrarily close to the identity transformation and the diffeomorphism $S$ arbitrarily close to $R$.
\end{Example}

A somewhat different example of this kind (with another form of the functions $\chi$, $\eta$, $P$ and without the observation that the diffeomorphism~$S$ can be made arbitrarily close to the rotation~$R$) was presented in the note~\cite{Sevryuk18}. To obtain an example for an open annulus, it is sufficient to disregard the boundary circles $\{\rho=1\}$ and $\{\rho=5\}$ in Example~\ref{exampletwo}.

Most probably, pairs of diffeomorphisms $R,S\colon\mM^2\to\mM^2$ with the intersection property such that the composition $RS$ also possesses the intersection property are an exception in some sense, but this issue is beyond the scope of this paper. Here we only note that the opposite situation where the composition $RS$ does not possess the intersection property (Example~\ref{exampletwo} is devoted to) is obviously persistent under small perturbations of $R$ and $S$.

\section{Non-intersections in dimension four and higher}\label{multi}

Now, for any $n$, let $L=\mT^n=\mR^n/2\pi\mZ^n$ be the standard $n$-torus, so that $T^\ast L\cong\mR^n\times\mT^n$ is a $2n$-dimensional cylinder. For $n\geq 2$ one can no longer assert that for any exact symplecto\-mor\-phism~$A$ of~$T^\ast\mT^n$ and any $n$-torus $\Gamma\subset T^\ast\mT^n$ homotopic to $\{0\}\times\mT^n$ ($0\in\mR^n$, $\Gamma\cong\mT^n$), the $n$-torus~$A\Gamma$ intersects $\Gamma$. Before giving appropriate counterexamples, we briefly describe a~general context of the problem.

One of the aspects of V.I.~Arnold's famous conjectures that initiated the development of symplectic topology (they go back to 1965 \cite{Arnold65}) concerns intersections of Hamiltonianly isotopic (or Hamiltonian homologous) Lagrangian submanifolds of a symplectic manifold. Recall that a submanifold $L$ of a symplectic manifold $\big(M^{2n},\omega^2\big)$ of dimension $2n$ is said to be \emph{Lagrangian} if $\dim L=n$ and $\omega^2|_L=0$. In the particular case where $M$ is the cotangent bundle of an $n$-dimensional manifold $\big(\omega^2=\rmd\lambda^1\big)$, a Lagrangian submanifold $L\subset M$ is said to be \emph{exact} if the $1$-form $\lambda^1|_L$ is not only closed but even exact. Two subsets of a symplectic manifold are said to be \emph{Hamiltonianly isotopic} if one of them can be obtained from the other by the phase flow mapping $\phi_0^1$ of a Hamiltonian (not necessarily autonomous) vector field~$V_t$. The general principle states that if two Lagrangian submanifolds $L_0$, $L_1$ of a symplectic manifold $\big(M,\omega^2\big)$ are Hamiltonianly isotopic then under certain additional (and not very restrictive) conditions they intersect: $L_0\cap L_1\neq\varnothing$, and the number of intersection points of~$L_0$ and~$L_1$ can often be estimated from below by various topological characteristics of $L_0$.

In particular, the presence of intersection points of $L_0$ and $L_1$ was established in the following cases (we will not dwell here on estimates for the number of intersection points):
\begin{enumerate}\itemsep=0pt
\item[a)] $M$ and $L_0$ are compact and the second relative homotopy group $\pi_2(M,L_0)$ vanishes \cite{Floer88};

\item[b)] $M=T^\ast L$ ($L$ being compact) and $L_0\cong L$ is the zero section of $T^\ast L$ \cite{Chaperon83,Hofer85,LaudenbachSikorav85} (in the paper \cite{Chaperon83}, only the case where $L=\mT^n$ was treated);

\item[c)] $L_0$ is compact and the symplectomorphism $\phi_0^1$ is sufficiently close to the identity transformation; in fact, case~c) follows from case~b) because a neighborhood of a Lagrangian submanifold $L\subset M$ is always symplectomorphic to a neighborhood of the zero section of~$T^\ast L$~\cite{Weinstein79};

\item[d)] $M=T^\ast L$ ($L$ being compact) and $L_0$ is compact and exact \cite{Gromov85}; by the way, under these conditions $L_0$ always intersects the zero section of $T^\ast L$ \cite{Gromov85}.
\end{enumerate}

All these results are discussed and the relevant references are given in, e.g., the works \cite{AbbondandoloSchlenk19,Arnold86,Audin14,Biran02,McDuffSalamon17}. In \cite{AbbondandoloSchlenk19,Biran02,McDuffSalamon17}, the subject is placed in the general context of the Floer homology theory which is very essential for the problems in question---according to \cite{Biran06}, ``whenever the principle of Lagrangian intersections fails (in the sense that a Lagrangian can be Hamiltonianly displaced), we obtain restrictions on the topology of the Lagrangian via computations in Floer homology''. Moreover, if $L_0$ and $L_1$ are two arbitrary compact and exact Lagrangian submanifolds of the cotangent bundle $T^\ast L$ of a compact manifold $L$ then the intersection of $L_0$ and $L_1$ is never empty (the condition of being Hamiltonianly isotopic is not needed), see the preprint \cite{Shelukhin19} and references therein, in particular, the paper~\cite{AbouzaidKragh18}.

Recall that a submanifold $L$ of a symplectic manifold $\big(M^{2n},\omega^2\big)$ is said to be \emph{coisotropic} if at any point $K\in L$, the tangent space $T_KL$ contains its skew-orthogonal complement (and, consequently, $\dim L\geq n$). Lagrangian submanifolds are just coisotropic submanifolds of the minimal possible dimension $n$. The parallel theory of unavoidable intersections of coisotropic submanifolds was founded by V.L.~Ginzburg~\cite{Ginzburg07}.

For non-Lagrangian submanifolds $L\subset M$ of a symplectic manifold $\big(M^{2n},\omega^2\big)$, some opposite statements are known.
\begin{enumerate}\itemsep=0pt
\item[A)] Let $L$ be compact, connected, of dimension $n$, and non-Lagrangian. Let also the normal fiber bundle of $L$ in $M$ admit a non-vanishing section. Then there exists a Hamiltonian vector field $V$ in $M$ nowhere tangent to $L$ (as one says, $L$ is \emph{infinitesimally displaceable})~\cite{LaudenbachSikorav94} (so that the image of $L$ under the phase flow mapping $\phi^\epsilon$ of $V$ does not intersect $L$ for any sufficiently small $\epsilon\neq 0$). The same conclusion was obtained earlier (but published slightly later) under more restrictive assumptions, namely, that $L$ is compact, connected, of dimension $n$, and non-Lagrangian, the Euler characteristic of $L$ vanishes, and there exists a field of Lagrangian subspaces along $L$ that is transversal to~$L$~\cite{Polterovich95}.

\item[B)] Let $L$ be compact, connected, and \emph{nowhere} coisotropic. Let also the normal fiber bundle of~$L$ in~$M$ admit a non-vanishing section. Then there exists a Hamiltonian vector field~$V$ in~$M$ nowhere tangent to~$L$~\cite{Gurel08}. In contrast with the non-Lagrangian case, the requirement here that $L$ is nowhere coisotropic cannot be relaxed to the condition that $L$ is (somewhere) non-coisotropic~\cite{Gurel08}.
\end{enumerate}

Now return to the $2n$-dimensional cylinder $T^\ast\mT^n$ with coordinates $\big(p\in\mR^n, q\in\mT^n\big)$, where~$p$ are the momenta corresponding to $q$. Given $v\in\mR^n$, the Hamiltonian vector field with the Hamilton function $v\cdot p$ is $v\partial_q$, and the phase flow mapping $\phi^t$ of $v\partial_q$ is the rotation \mbox{$(p,q)\mapsto(p,q+tv)$} by the ``multidimensional'' $p$-independent angle $tv$. Such rotations are simplest instances of exact symplectomorphisms of $T^\ast\mT^n$.

\begin{Example}\label{examplemulti}
A non-Lagrangian $n$-torus $L_0\subset T^\ast\mT^n$ \textup($n\geq 2$\textup) arbitrarily close to the torus $\big\{p=p^0\big\}$, $p^0\in\mR^n$, and a rotation $R\colon (p,q)\mapsto(p,q+v)$, $v\in\mR^n$, such that $L_0\cap\big\{p=p^0\big\}=\varnothing$ and $L_0\cap RL_0=\varnothing$.

Let $L_0=\big\{p=p^0+\varepsilon f(q_1)\big\}$ where $\varepsilon\neq 0$ is arbitrarily small, $f_1(q_1)=\sin^kq_1$, $f_2(q_1)=\cos q_1$, $f_i(q_1)\equiv 0$ for $3\leq i\leq n$, and $k$ is equal to $1$ or $3$. The $n$-torus $L_0$ is not Lagrangian since $(\rmd p\wedge\rmd q)|_{L_0} = -\varepsilon(\sin q_1) \, \rmd q_1\wedge\rmd q_2$ (although one cannot say that $L_0$ is nowhere Lagrangian), it does not intersect the torus $\big\{p=p^0\big\}$ (in particular, if $p^0=0$ then $L_0$ does not intersect the zero section $\{p=0\}$ of $T^\ast\mT^n$), and $RL_0=\big\{p=p^0+\varepsilon f(q_1-v_1)\big\}$, so that $RL_0$ intersects $L_0$ if and only if $v_1\equiv 0\bmod 2\pi$.

In this setup, there is however the following significant difference between the cases $k=1$ and $k=3$. The tangent space to the torus $L_0$ at a point $\bigl(p^0+\varepsilon f(q_1), q\bigr)$ is spanned by the vectors $\varepsilon f'_1(q_1)\partial_{p_1}+\varepsilon f'_2(q_1)\partial_{p_2}+\partial_{q_1}$, $\partial_{q_2}$, $\ldots$, $\partial_{q_n}$. Thus, if $k=1$ then the vector field $v\partial_q$ with $v_1\neq 0$ is nowhere tangent to $L_0$. On the other hand, if $k=3$ then the vector field $v\partial_q$ with $v_1\neq 0$ is tangent to $L_0$ at the points where $\sin q_1=0$, i.e., along the two $(n-1)$-tori corresponding to $q_1=0$ and $q_1=\pi$. These are precisely the $(n-1)$-tori along which the $2$-form $(\rmd p\wedge\rmd q)|_{L_0}$ vanishes.

Of course, according to the theorem by F.~Laudenbach and J.-C.~Sikorav \cite{LaudenbachSikorav94} mentioned above, a Hamiltonian vector field in $T^\ast\mT^n$ nowhere tangent to $L_0$ exists in the case $k=3$ as well. For instance, let $k=3$ and consider the Hamiltonian vector field $V = a\partial_{q_1}+b(\cos q_1)\partial_{p_1}$ with the Hamilton function $ap_1-b\sin q_1$ ($a,b\in\mR$). If $ab\neq 0$ then $V$ is nowhere tangent to~$L_0$. Indeed, if $V$ is tangent to $L_0$ at a point $\bigl(p^0+\varepsilon f(q_1), q\bigr)$ then $f'_2(q_1) = -\sin q_1 = 0$ and $\varepsilon af'_1(q_1) = 3\varepsilon a\sin^2q_1\cos q_1 = b\cos q_1$ but these two equalities cannot hold simultaneously.
\end{Example}

We do not consider here multidimensional generalizations of the intersection property. Recent works on this subject are exemplified by the papers~\cite{LiuXing22,YangLi20}.

\section{Topologically non-equivalent phase portraits}\label{equivalence}

While the two preceding sections were devoted to intersection patterns of diffeomorphisms of multidimensional annuli, in this section we deal with a completely different range of problems, namely, with the topological classification of the phase portraits of planar vector fields. Let $(x,y)$ be Cartesian coordinates in~$\mR^2$. We will be interested in one-parameter families of analytic autonomous ordinary differential equations of the form
\begin{equation}
\dot{x}=f(x,y), \qquad \dot{y}=\mu g(x,y)\label{mu}
\end{equation}
($\mu>0$) in the closed upper half-plane $\{(x,y) \mid y\geq 0\}$ such that the topological type of the phase portrait of the equation~\eqref{mu} changes as the parameter $\mu$ is varied from $0$ to $+\infty$.

It is very easy to present an example of such a family in the entire plane $\mR^2$ or in the open upper half-plane $\{(x,y) \mid y>0\}$. In the case of the entire plane $\mR^2$, one may consider, e.g., the family of linear equations
\begin{equation}
\dot{x}=3(x+y), \qquad \dot{y}=-\mu(4x+3y)
\label{mymu}
\end{equation}
(cf.\ \cite{AVK22}). The characteristic polynomial of the matrix of the system~\eqref{mymu} is $\lambda^2+3(\mu-1)\lambda+3\mu$, and the eigenvalues of this matrix are
\[
\frac{3(1-\mu)\pm[D(\mu)]^{1/2}}{2}, \qquad D(\mu) = 9(\mu-1)^2-12\mu = 3(3\mu-1)(\mu-3).
\]
Thus, the equilibrium $\{x=y=0\}$ of the system~\eqref{mymu} is an unstable node for $0<\mu<1/3$, an unstable degenerate (Jordan) node for $\mu=1/3$ (with eigenvalues $1$, $1$), an unstable focus for $1/3<\mu<1$, a center for $\mu=1$ \big(with eigenvalues $\pm 3^{1/2}\rmi$\big), a stable focus for $1<\mu<3$, a stable degenerate (Jordan) node for $\mu=3$ (with eigenvalues $-3$, $-3$), and a stable node for $\mu>3$. Consequently, the topological types of the phase portrait of the system~\eqref{mymu} for $0<\mu<1$, $\mu=1$, and $\mu>1$ are different. In the case of the open upper half-plane $\{(x,Y) \mid Y>0\}$, one may perform the coordinate change $Y=\rme^y$ where the variables $(x,y)\in\mR^2$ satisfy the system~\eqref{mymu}, cf.\ \cite{AVK22}:
\[
\dot{x}=3(x+\ln Y), \qquad \dot{Y}=-\mu Y(4x+3\ln Y).
\]

\begin{Example}\label{exampleportraits1}
A one-parameter family of analytic autonomous ordinary differential equations of the form~\eqref{mu} in the closed upper half-plane $\{(x,y) \mid y\geq 0\}$ such that the corresponding phase portraits for $0<\mu<1$ and for $\mu\geq 1$ are topologically non-equivalent but in a neighborhood of the origin, the phase portraits for all values $\mu>0$ of the parameter are topologically equivalent.

Consider the family
\begin{equation}
\dot{x}=x^2-y^2, \qquad \dot{y}=\mu xy.
\label{mu1}
\end{equation}
For any $\mu$, the origin $O = \{x=y=0\}$ (an equilibrium) and the rays $\{y=0, x<0\}$, $\{y=0, x>0\}$ are phase curves of the system~\eqref{mu1}. The phase curves in the open upper half-plane $\{(x,y) \mid y>0\}$ can be described as follows.

Let $a\in\mR$. It is easy to verify that the curves
\[
\fI_{\mu,a} = \left\{ (x,y) \Biggm| x^2 = ay^{2/\mu}+\frac{y^2}{1-\mu} \right\}
\]
are invariant under the phase flow of the system~\eqref{mu1} for any positive $\mu\neq 1$ and any $a\in\mR$ (if $\mu>1$ then it is supposed that $a>0$) and that the curves
\[
\fJ_a = \bigl\{ (x,y) \bigm| x^2 = y^2(a-2\ln y) \bigr\}
\]
(where $y^2\ln y|_{y=0}=0$) are invariant under the phase flow of the system~\eqref{mu1} for $\mu=1$ and any $a\in\mR$. Thus, the phase curves of the system~\eqref{mu1} in the domain $\{y>0\}$ are as follows.

For $0<\mu<1$ (Figure~\ref{fig41}):
\begin{itemize}\itemsep=0pt
\item[--] the curves $\fI_{\mu,a}\setminus O$ for $a<0$; note that each curve $\fI_{\mu,a}$ for $0<\mu<1$ and $a<0$ is closed and intersects the ordinate axis $\{x=0\}$ at two points with the ordinates $y=0$ and $y=[a(\mu-1)]^{\mu/(2\mu-2)}$, so that the phase curve $\fI_{\mu,a}\setminus O$ is homoclinic to the equilibrium $O$,
\item[--] the curves $\fI_{\mu,a} \cap \{x<0\}$ and $\fI_{\mu,a} \cap \{x>0\}$ for $a\geq 0$; note that the curves $\fI_{\mu,a} \cap \{x<0\}$ come from infinity and are forward asymptotic to $O$ while the curves $\fI_{\mu,a} \cap \{x>0\}$ are backward asymptotic to $O$ and go to infinity; for $a=0$ these curves are straight rays.
\end{itemize}

\begin{figure}[t!]\centering
\includegraphics{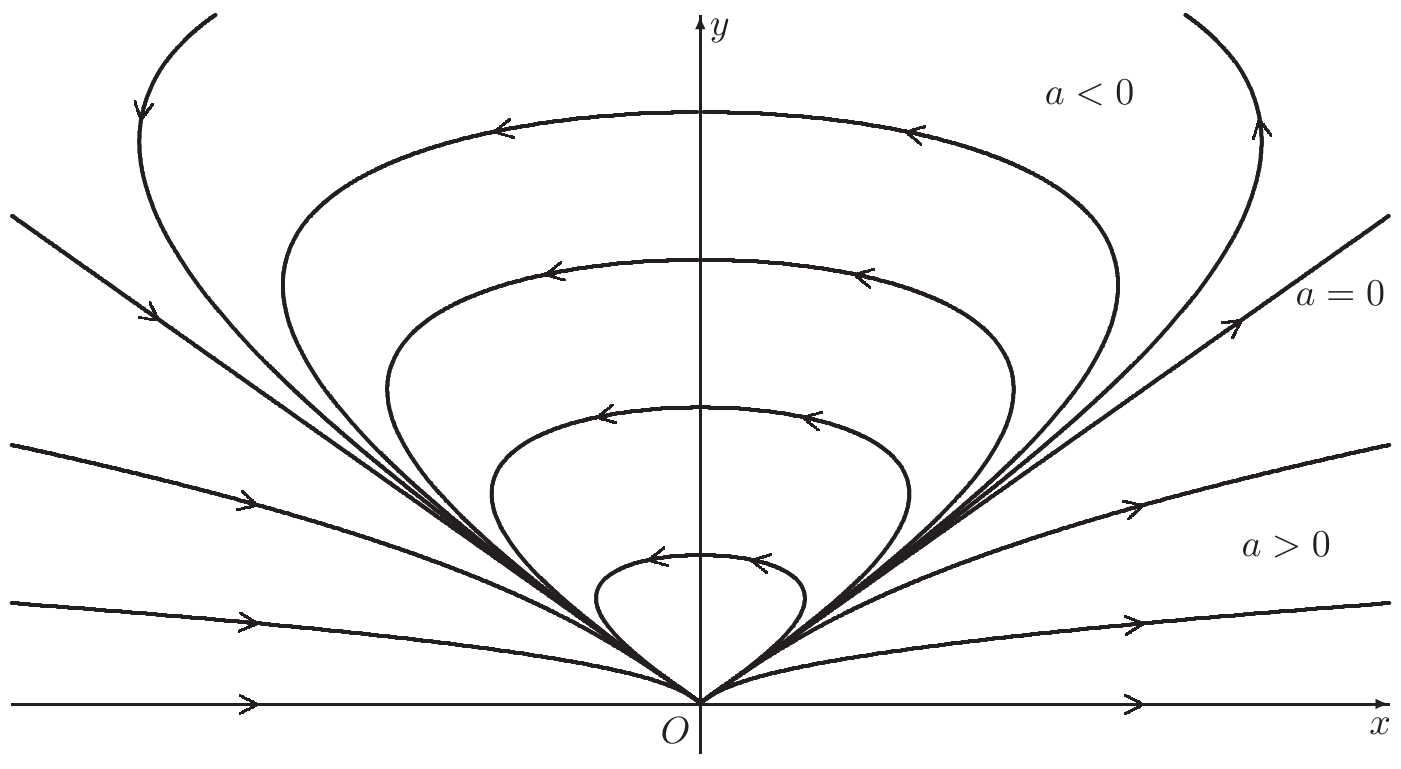}
\caption{The phase portrait of the system~\eqref{mu1} for $\mu=0.5$ with $\fI_{1/2,a} = \{x^2=ay^4+2y^2\}$.}\label{fig41}
\end{figure}

For $\mu=1$ (Figure~\ref{fig42}):
\begin{itemize}\itemsep=0pt
\item[--] the curves $\fJ_a\setminus O$ for $a\in\mR$; note that each curve $\fJ_a$ is closed and intersects the ordinate axis $\{x=0\}$ at two points with the ordinates $y=0$ and $y=\rme^{a/2}$, so that the phase curve $\fJ_a\setminus O$ is homoclinic to the equilibrium $O$.
\end{itemize}

\begin{figure}[th!]\centering
\includegraphics{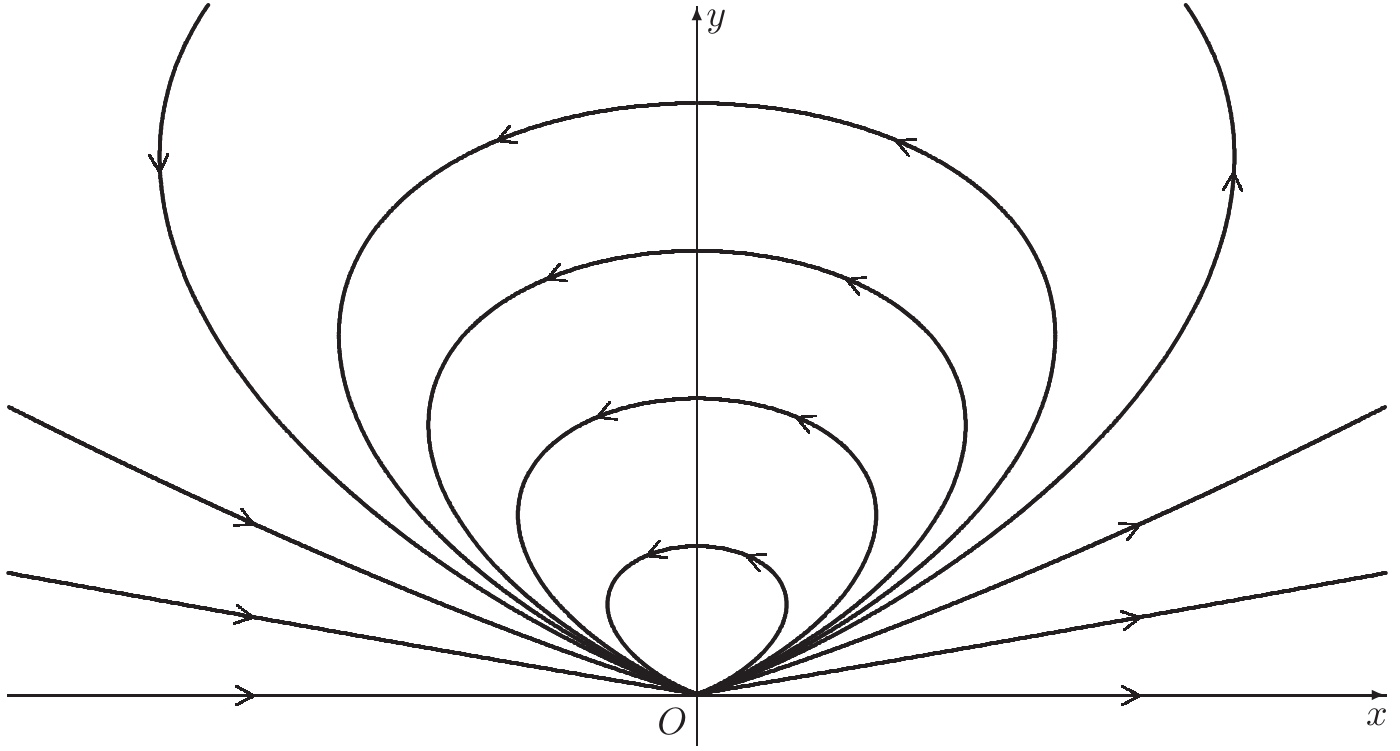}
\caption{The phase portrait of the system~\eqref{mu1} for $\mu=1$.}\label{fig42}
\end{figure}

For $\mu>1$ (Figure~\ref{fig43}):
\begin{itemize}\itemsep=0pt
\item[--] the curves $\fI_{\mu,a}\setminus O$ for $a>0$; note that each curve $\fI_{\mu,a}$ for $\mu>1$ and $a>0$ is closed and intersects the ordinate axis $\{x=0\}$ at two points with the ordinates $y=0$ and $y=[a(\mu-1)]^{\mu/(2\mu-2)}$, so that the phase curve $\fI_{\mu,a}\setminus O$ is homoclinic to the equilibrium $O$.
\end{itemize}

\begin{figure}[th!]\centering
\includegraphics{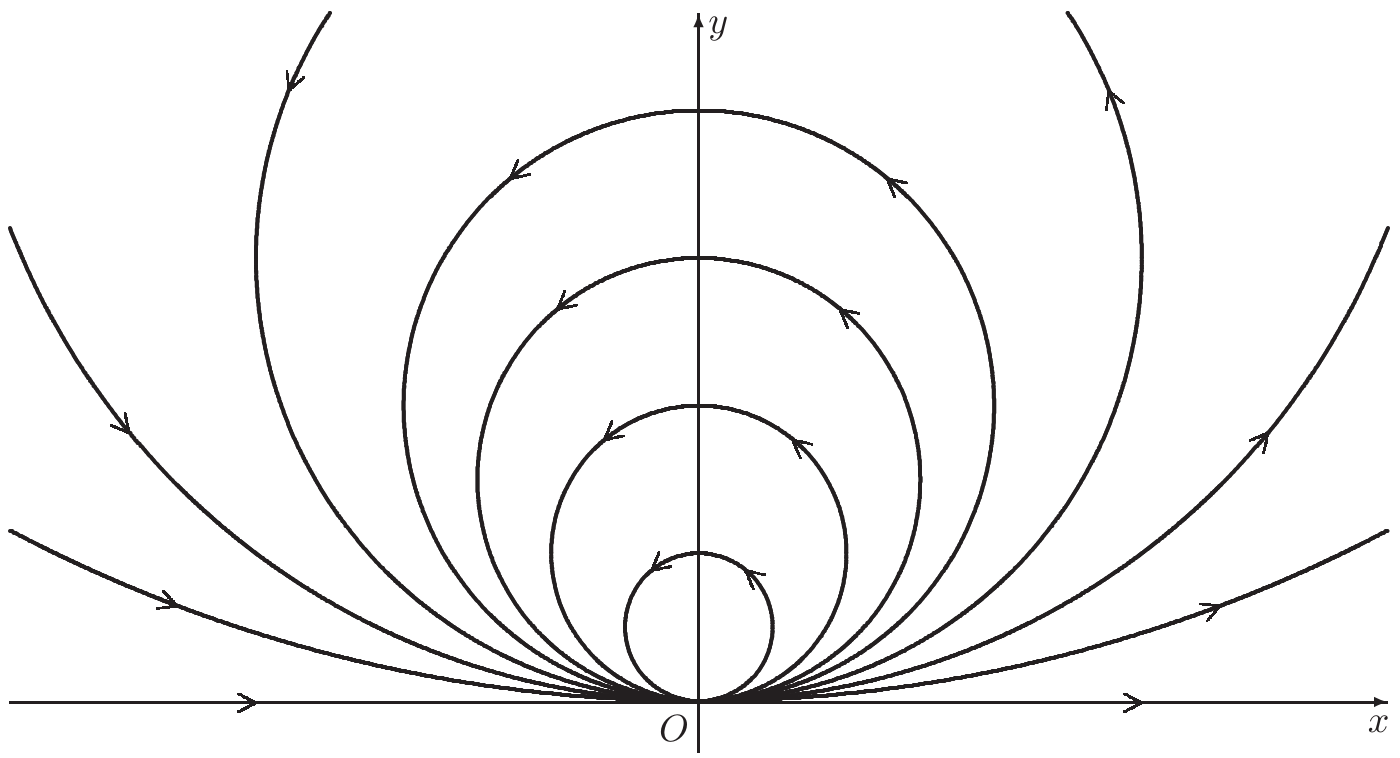}
\caption{The phase portrait of the system~\eqref{mu1} for $\mu=2$. The curves $\fI_{2,a}$ are the circles $x^2=ay-y^2$.}\label{fig43}
\end{figure}

It is clear that the phase portraits of the system~\eqref{mu1} for $0<\mu<1$ and for $\mu\geq 1$ are globally topologically non-equivalent. Indeed, for $0<\mu<1$, there are infinitely many phase curves joining the origin with infinity whereas for $\mu>1$, such phase curves in the domain $\{y>0\}$ are absent. However, the phase portraits for all $\mu>0$ are topologically equivalent locally, i.e., in a neighborhood of the origin.
\end{Example}

\begin{Example}\label{exampleportraits2}
A one-parameter family of analytic autonomous ordinary differential equations of the form~\eqref{mu} in the closed upper half-plane $\{(x,y) \mid y\geq 0\}$ such that the corresponding phase portraits for $0<\mu\leq 1$ and for $\mu>1$ are topologically non-equivalent both globally (in the entire closed upper half-plane) and locally (in a neighborhood of the origin).

Consider the family
\begin{equation}
\dot{x}=x^2+y^2, \qquad \dot{y}=\mu xy.\label{mu2}
\end{equation}
As in the case of the family~\eqref{mu1}, for any $\mu$, the origin $O = \{x=y=0\}$ (an equilibrium) and the rays $\{y=0, x<0\}$, $\{y=0, x>0\}$ are phase curves of the system~\eqref{mu2}. The phase curves in the open upper half-plane $\{(x,y) \mid y>0\}$ can be described as follows.

Let $a\in\mR$. It is easy to verify that the curves
\[
\fK_{\mu,a} = \left\{ (x,y) \Biggm| x^2 = ay^{2/\mu}+\frac{y^2}{\mu-1} \right\}
\]
are invariant under the phase flow of the system~\eqref{mu2} for any positive $\mu\neq 1$ and any $a\in\mR$ (if $0<\mu<1$ then it is supposed that $a>0$) and that the curves
\[
\fL_a = \bigl\{ (x,y) \bigm| x^2 = y^2(a+2\ln y) \bigr\}
\]
(where $y^2\ln y|_{y=0}=0$) are invariant under the phase flow of the system~\eqref{mu2} for $\mu=1$ and any $a\in\mR$. Thus, the phase curves of the system~\eqref{mu2} in the domain $\{y>0\}$ are as follows.

For $0<\mu<1$ (Figure~\ref{fig44}):
\begin{itemize}\itemsep=0pt
\item[--] the curves $\fK_{\mu,a}$ for $a>0$; note that each phase curve $\fK_{\mu,a}$ for $0<\mu<1$ and $a>0$ comes from infinity, goes to infinity, and intersects the ordinate axis $\{x=0\}$ at the only point with the ordinate $y=[a(1-\mu)]^{\mu/(2\mu-2)}$.
\end{itemize}

\begin{figure}[t!]\centering
\includegraphics{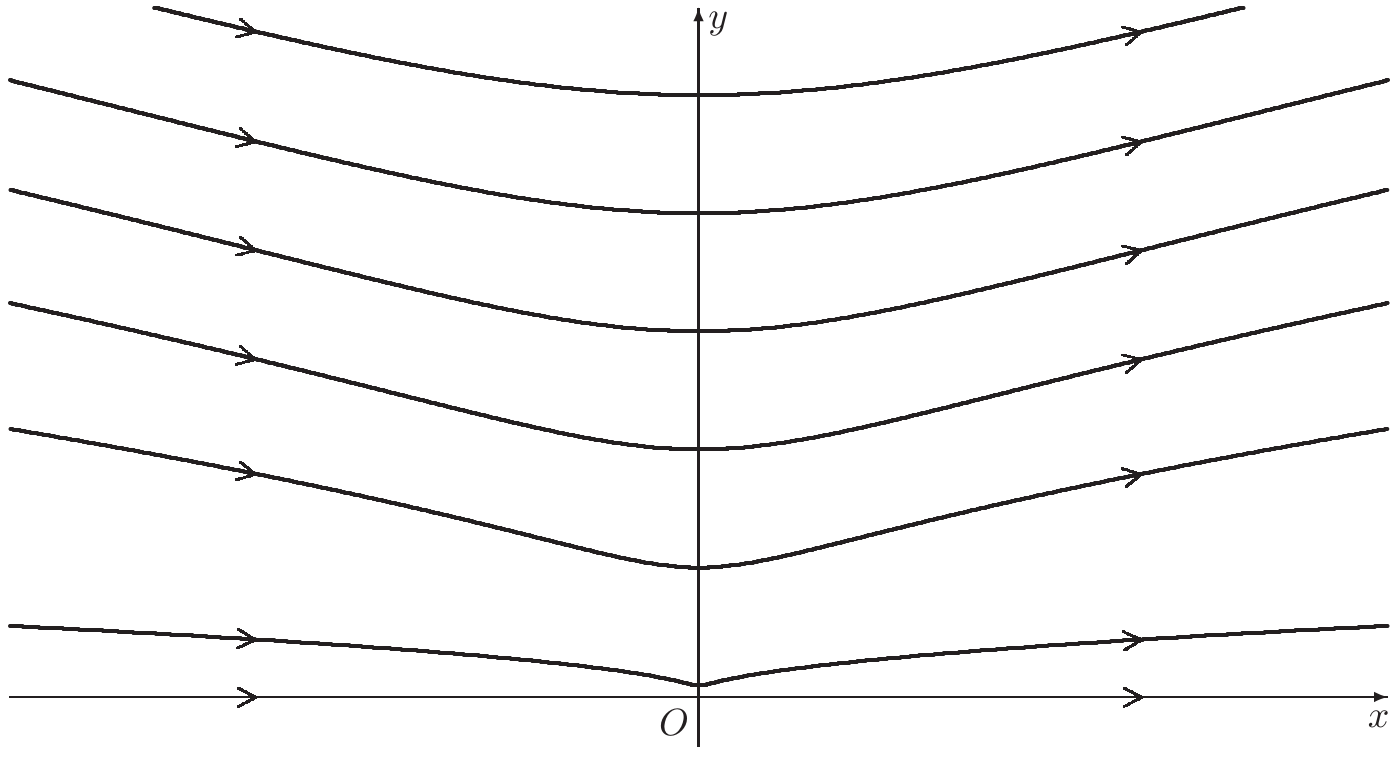}
\caption{The phase portrait of the system~\eqref{mu2} for $\mu=0.5$ with $\fK_{1/2,a} = \{x^2=ay^4-2y^2\}$.}\label{fig44}
\end{figure}

For $\mu=1$ (Figure~\ref{fig45}):
\begin{itemize}\itemsep=0pt
\item[--] the curves $\fL_a$ for $a\in\mR$; note that each phase curve $\fL_a$ for $\mu=1$ comes from infinity, goes to infinity, and intersects the ordinate axis $\{x=0\}$ at the only point with the ordinate $y=\rme^{-a/2}$.
\end{itemize}

\begin{figure}[th!]\centering
\includegraphics{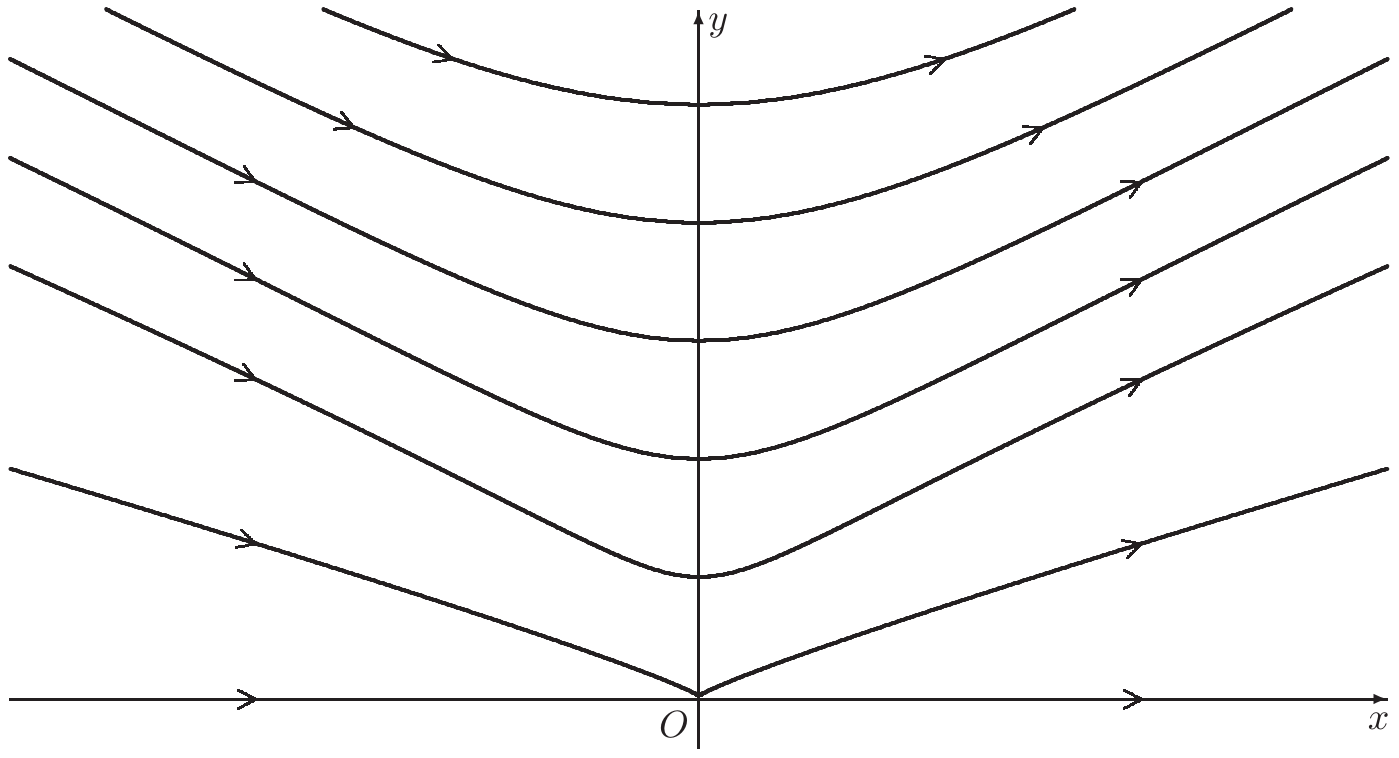}
\caption{The phase portrait of the system~\eqref{mu2} for $\mu=1$.}\label{fig45}
\end{figure}

For $\mu>1$ (Figure~\ref{fig46}):
\begin{itemize}\itemsep=0pt
\item[--] the curves $\fK_{\mu,a}$ for $a<0$; note that each phase curve $\fK_{\mu,a}$ for $\mu>1$ and $a<0$ comes from infinity, goes to infinity, and intersects the ordinate axis $\{x=0\}$ at the only point with the ordinate $y=[a(1-\mu)]^{\mu/(2\mu-2)}$,
\item[--] the curves $\fK_{\mu,a} \cap \{x<0\}$ and $\fK_{\mu,a} \cap \{x>0\}$ for $a\geq 0$; note that the curves $\fK_{\mu,a} \cap \{x<0\}$ come from infinity and are forward asymptotic to~$O$ while the curves $\fK_{\mu,a} \cap \{x>0\}$ are backward asymptotic to $O$ and go to infinity; for $a=0$ these curves are straight rays.
\end{itemize}

\begin{figure}[th!]\centering
\includegraphics{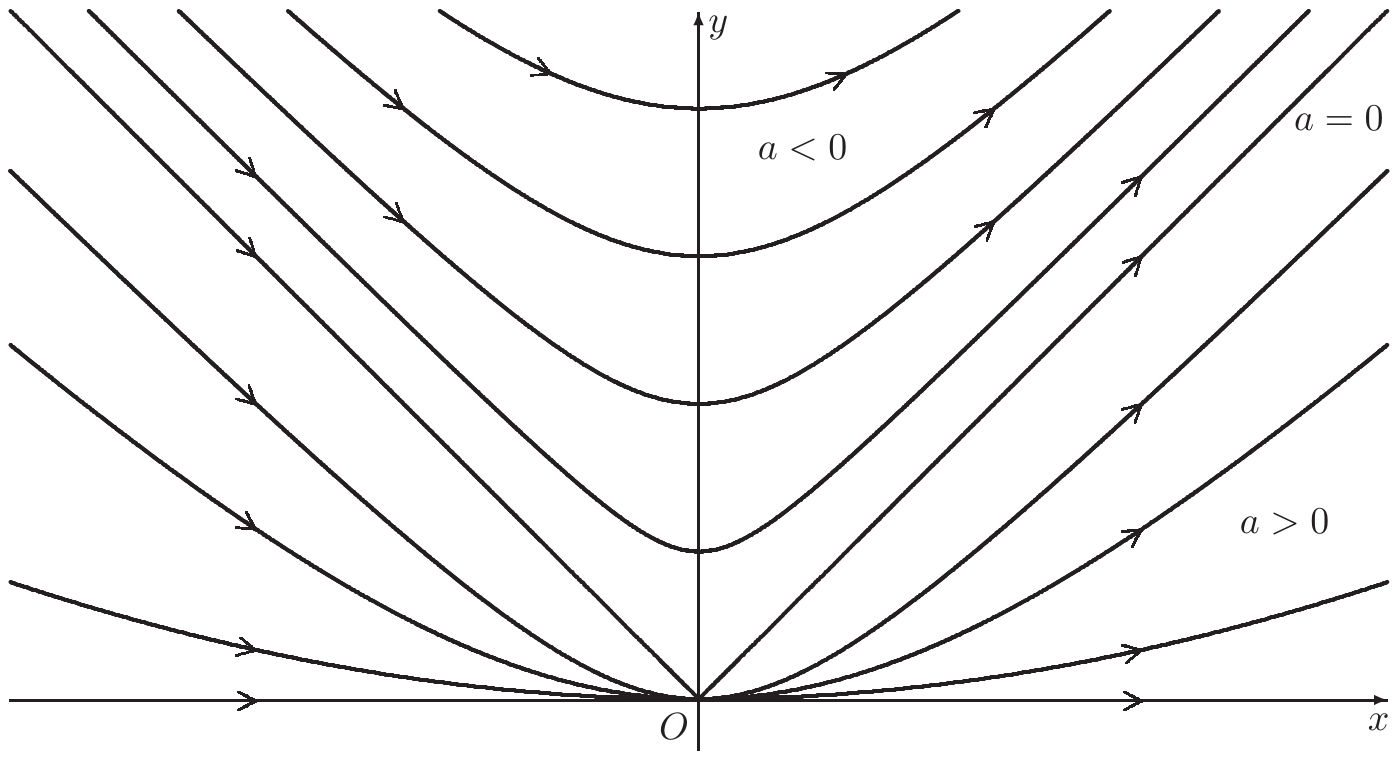}
\caption{The phase portrait of the system~\eqref{mu2} for $\mu=2$. The curves $\fK_{2,a}$ for $a\neq 0$ are the equilateral hyperbolas $x^2=ay+y^2$.}\label{fig46}
\end{figure}

It is clear that the phase portraits of the system~\eqref{mu2} for $0<\mu\leq 1$ and for $\mu>1$ are topologically equivalent neither globally nor locally. Indeed, for $\mu>1$, there are infinitely many phase curves forward asymptotic to $O$ and infinitely many phase curves backward asymptotic to $O$ whereas for $0<\mu\leq 1$, such phase curves in the domain $\{y>0\}$ are absent.
\end{Example}

Of course, Examples~\ref{exampleportraits1} and~\ref{exampleportraits2} can be seen as two variations of the same example. Note that for all $\mu$, the systems~\eqref{mu1} and~\eqref{mu2} are reversible with respect to the involution $(x,y)\mapsto(-x,y)$. Note also that since the right-hand sides of the systems~\eqref{mu1} and~\eqref{mu2} are quadratic forms, the phase portraits of these systems are invariant under the homotheties $H_k\colon (x,y)\mapsto(kx,ky)$ ($k>0$). Indeed, $H_k\fI_{\mu,a} = \fI_{\mu,\kappa a}$, $H_k\fK_{\mu,a} = \fK_{\mu,\kappa a}$ with $\kappa=k^{(2\mu-2)/\mu}$ and $H_k\fJ_a = \fJ_{a+2\ln k}$, $H_k\fL_a = \fL_{a-2\ln k}$.

One more remark is that the vector fields in both Examples~\ref{exampleportraits1} and~\ref{exampleportraits2} are incomplete. To see this, it suffices to notice that in both the examples, the induced dynamics in the invariant line $\{y=0\}$ is described by the equation $\dot{x}=x^2$ ($x^2\partial_x$ is a ``prototype'' of incomplete vector fields in the real line $\mR$). To obtain complete examples, it is sufficient to divide the right-hand sides of the systems~\eqref{mu1} and~\eqref{mu2} by, e.g., $1+x^2+y^2$. This would not alter the phase portraits, but we would get equations~\eqref{mu} with $|f(x,y)|\leq 1$ and $|g(x,y)|\leq 1/2$ for any $x$ and $y$. The corresponding vector fields would obviously be complete.

In fact, the difficulty in searching for families of systems like those in Examples~\ref{exampleportraits1} and~\ref{exampleportraits2} comes from two requirements: that the phase space has a boundary and the vector fields are tangent to it and that the vector fields are analytic. Without these restrictions, examples would be straightforward to construct (see the discussions at the beginning of this section and in the post \cite{AVK22} where a $C^\infty$ example is presented).

\subsection*{Acknowledgements}

I am grateful to L.V.~Polterovich for interesting discussions of some aspects of symplectic topology and a number of useful references (in particular, the preprint~\cite{Shelukhin19}). Special thanks go also to V.M.~Zyuzkov. I am appreciative for helpful critical remarks of anonymous referees and for the editor's comments.

\pdfbookmark[1]{References}{ref}
\LastPageEnding


\begin{thebibliography}{99}
\footnotesize\itemsep=0pt

\bibitem{AbbondandoloSchlenk19}
Abbondandolo A., Schlenk F., Floer homologies, with applications,
 \href{https://doi.org/10.1365/s13291-018-0193-x}{\textit{Jahresber. Dtsch. Math.-Ver.}} \textbf{121} (2019), 155--238,
 \href{https://arxiv.org/abs/1709.00297}{arXiv:1709.00297}.

\bibitem{AbouzaidKragh18}
Abouzaid M., Kragh T., Simple homotopy equivalence of nearby {L}agrangians,
 \href{https://doi.org/10.4310/ACTA.2018.v220.n2.a1}{\textit{Acta Math.}} \textbf{220} (2018), 207--237, \href{https://arxiv.org/abs/1603.05431}{arXiv:1603.05431}.

\bibitem{Arnold65}
Arnold V.I., Sur une propri\'et\'e topologique des applications globalement
 canoniques de la m\'ecanique classique, \textit{C.~R.~Acad. Sci. Paris}
 \textbf{261} (1965), 3719--3722, {E}nglish translation: On a topological
 property of globally canonical maps in classical mechanics, in Vladimir
 I.~Arnold. Collected Works, Vol.~1, Editors A.B.~Givental, B.A.~Khesin,
 J.E.~Marsden, A.N.~Varchenko, V.A.~Vassiliev, O.Ya.~Viro, V.M.~Zakalyukin,
 \href{https://doi.org/10.1007/978-3-642-01742-1}{Springer-Verlag}, Berlin,
 2009, 481--485.

\bibitem{Arnold86}
Arnold V.I., First steps of symplectic topology, \href{https://doi.org/10.1070/RM1986v041n06ABEH004221}{\textit{Russian Math. Surveys}}
 \textbf{41} (1986), no.~6, 1--21.

\bibitem{Arnold96}
Arnold V.I., Topological problems of the theory of wave propagation,
 \href{https://doi.org/10.1070/RM1996v051n01ABEH002734}{\textit{Russian Math. Surveys}} \textbf{51} (1996), no.~1, 1--47.

\bibitem{Arnold14}
Arnold V.I., About {V}ladimir {A}bramovich {R}okhlin, in Arnold: Swimming
 Against the Tide, \textit{AMS Non-Series Monographs}, Vol.~86, Editors B.A.~Khesin, S.L.~Tabachnikov, \href{https://doi.org/10.1090/mbk/086}{Amer. Math. Soc.}, Providence, RI, 2014, 67--75.

\bibitem{Audin14}
Audin M., Vladimir {I}gorevich {A}rnold and the invention of symplectic
 topology, in Contact and Symplectic Topology, \textit{Bolyai Soc. Math.
 Stud.}, Vol.~26, \href{https://doi.org/10.1007/978-3-319-02036-5_1}{J\'anos Bolyai Math. Soc.}, Budapest, 2014, 1--25.

\bibitem{Biran02}
Biran P., Geometry of symplectic intersections, in Proceedings of the
 {I}nternational {C}ongress of {M}athematicians, {V}ol.~{II} ({B}eijing,
 2002), \href{https://doi.org/10.1142/4962}{Higher Ed. Press}, Beijing, 2002, 241--255, \href{https://arxiv.org/abs/math.SG/0304260}{arXiv:math.SG/0304260}.

\bibitem{Biran06}
Biran P., Lagrangian non-intersections, \href{https://doi.org/10.1007/s00039-006-0560-0}{\textit{Geom. Funct. Anal.}} \textbf{16}
 (2006), 279--326, \href{https://arxiv.org/abs/math.SG/0412110}{arXiv:math.SG/0412110}.

\bibitem{Boss20}
Boss V., Counterexamples and paradoxes, 3rd ed., \textit{Lectures on Math.}, Vol.~12, Librokom, Moscow, 2020.

\bibitem{Chaperon83}
Chaperon M., Quelques questions de g\'eom\'etrie symplectique (d'apr\`es, entre
 autres, {P}oincar\'e, {A}rnold, {C}onley et {Z}ehnder), \textit{Ast\'erisque}
 \textbf{105} (1983), 231--249.

\bibitem{Floer88}
Floer A., Morse theory for {L}agrangian intersections, \href{https://doi.org/10.4310/jdg/1214442477}{\textit{J.~Differential
 Geom.}} \textbf{28} (1988), 513--547.

\bibitem{GelbaumOlmsted03}
Gelbaum B.R., Olmsted J.M.H., Counterexamples in analysis, Dover Publications,
 Inc., Mineola, NY, 2003.

\bibitem{Ginzburg07}
Ginzburg V.L., Coisotropic intersections, \href{https://doi.org/10.1215/S0012-7094-07-14014-6}{\textit{Duke Math.~J.}} \textbf{140}
 (2007), 111--163, \href{https://arxiv.org/abs/math.SG/0605186}{arXiv:math.SG/0605186}.

\bibitem{Gromov85}
Gromov M., Pseudo holomorphic curves in symplectic manifolds, \href{https://doi.org/10.1007/BF01388806}{\textit{Invent.
 Math.}} \textbf{82} (1985), 307--347.

\bibitem{Gurel08}
G\"urel B.Z., Totally non-coisotropic displacement and its applications to
 {H}amiltonian dynamics, \href{https://doi.org/10.1142/S0219199708003198}{\textit{Commun. Contemp. Math.}} \textbf{10} (2008),
 1103--1128, \href{https://arxiv.org/abs/math.SG/0702091}{arXiv:math.SG/0702091}.

\bibitem{Hofer85}
Hofer H., Lagrangian embeddings and critical point theory, \href{https://doi.org/10.1016/S0294-1449(16)30394-8}{\textit{Ann. Inst.
 H.~Poincar\'e Anal. Non Lin\'eaire}} \textbf{2} (1985), 407--462.

\bibitem{HuangLiLiu18}
Huang P., Li X., Liu B., Invariant curves of smooth quasi-periodic mappings,
 \href{https://doi.org/10.3934/dcds.2018006}{\textit{Discrete Contin. Dyn. Syst.}} \textbf{38} (2018), 131--154,
 \href{https://arxiv.org/abs/1705.08762}{arXiv:1705.08762}.

\bibitem{HuangLiLiu22}
Huang P., Li X., Liu B., Invariant curves of almost periodic twist mappings,
 \href{https://doi.org/10.1007/s10884-021-10033-1}{\textit{J.~Dynam. Differential Equations}} \textbf{34} (2022), 1997--2033,
 \href{https://arxiv.org/abs/1606.08938}{arXiv:1606.08938}.

\bibitem{LaudenbachSikorav85}
Laudenbach F., Sikorav J.-C., Persistance d'intersection avec la section nulle
 au cours d'une isotopie hamiltonienne dans un fibr\'e cotangent,
 \href{https://doi.org/10.1007/BF01388807}{\textit{Invent. Math.}} \textbf{82} (1985), 349--357.

\bibitem{LaudenbachSikorav94}
Laudenbach F., Sikorav J.-C., Hamiltonian disjunction and limits of {L}agrangian
 submanifolds, \href{https://doi.org/10.1155/S1073792894000176}{\textit{Int. Math. Res. Not.}} \textbf{1994} (1994), 161~ff.,
 8~pages.

\bibitem{LiuXing22}
Liu C., Xing J., A new proof of {M}oser's theorem, \href{https://doi.org/10.11948/20220161}{\textit{J.~Appl. Anal.
 Comput.}} \textbf{12} (2022), 1679--1701.

\bibitem{MaXu22}
Ma Z., Xu J., A {KAM} theorem for quasi-periodic non-twist mappings and its
 application, \href{https://doi.org/10.3934/dcds.2022013}{\textit{Discrete Contin. Dyn. Syst.}} \textbf{42} (2022),
 3169--3185.

\bibitem{McDuffSalamon17}
McDuff D., Salamon D., Introduction to symplectic topology, 3rd ed., \textit{Oxf. Grad.
 Texts Math.}, \href{https://doi.org/10.1093/oso/9780198794899.001.0001}{Oxford University Press}, Oxford, 2017.

\bibitem{Moser62}
Moser J., On invariant curves of area-preserving mappings of an annulus,
 \textit{Nachr. Akad. Wiss. G\"{o}ttingen Math.-Phys. Kl.~II} \textbf{1962}
 (1962), 1--20.

\bibitem{Moser01}
Moser J., Remark on the paper ``{O}n invariant curves of area-preserving
 mappings of an annulus'', \href{https://doi.org/10.1070/RD2001v006n03ABEH000181}{\textit{Regul. Chaotic Dyn.}} \textbf{6} (2001),
 337--338.

\bibitem{Polterovich95}
Polterovich L., An obstacle to non-{L}agrangian intersections, in The {F}loer
 Memorial Volume, \textit{Progr. Math.}, Vol.~133, \href{https://doi.org/10.1007/978-3-0348-9217-9}{Birkh\"auser}, Basel, 1995,
 575--586.

\bibitem{Russmann83}
R\"ussmann H., On the existence of invariant curves of twist mappings of an
 annulus, in Geometric Dynamics ({R}io de {J}aneiro, 1981), \textit{Lecture
 Notes in Math.}, Vol.~1007, \href{https://doi.org/10.1007/BFb0061441}{Springer}, Berlin, 1983, 677--718.

\bibitem{Sevryuk18}
Sevryuk M.B., On the composition of mappings possessing the intersection
 property, \href{https://doi.org/10.20861/2312-8267-2018-45-003}{\textit{Sci. Tech. Educ.}} (2018), no.~4, 6--9.

\bibitem{Shafarevich05}
Shafarevich I.R., Basic notions of algebra, \textit{Encyclopaedia Math. Sci.},
 Vol.~11, \href{https://doi.org/10.1007/b137643}{Springer-Verlag}, Berlin, 2005.

\bibitem{Shelukhin19}
Shelukhin E., Symplectic cohomology and a conjecture of {V}iterbo,
 \textit{Geom. Funct. Anal.}, {t}o appear, \href{https://arxiv.org/abs/1904.06798}{arXiv:1904.06798}.

\bibitem{SiegelMoser95}
Siegel C.L., Moser J.K., Lectures on celestial mechanics, \textit{Classics Math.},
 \href{https://doi.org/10.1007/978-3-642-87284-6}{Springer-Verlag}, Berlin, 1995.

\bibitem{AVK22}
The post by the user {AVK} of {M}arch~5 (edited {M}arch~7), 2022, available at
 \url{https://mathoverflow.net/questions/417435/special-topological-equivalence}.

\bibitem{Weinstein79}
Weinstein A., Lectures on symplectic manifolds, \textit{CBMS Reg. Conf. Ser.
 Math.}, Vol.~29, \href{https://doi.org/10.1090/cbms/029}{Amer. Math. Soc.}, Providence, R.I., 1979.

\bibitem{XiaYu22}
Xia Z., Yu P., A fixed point theorem for twist maps, \href{https://doi.org/10.3934/dcds.2022045}{\textit{Discrete Contin.
 Dyn. Syst.}} \textbf{42} (2022), 4051--4059, \href{https://arxiv.org/abs/2106.06374}{arXiv:2106.06374}.

\bibitem{Yakovenko14}
Yakovenko S.Yu., Vladimir {I}gorevich {A}rnold: {A} view from the rear bench, in
 Arnold: Swimming Against the Tide, \textit{AMS Non-Series Monographs},
 Vol.~86, Editors B.A.~Khesin, S.L.~Tabachnikov, \href{https://doi.org/10.1090/mbk/086}{Amer. Math. Soc.}, Providence,
 RI, 2014, 197--203.

\bibitem{YangLi20}
Yang L., Li X., Existence of periodically invariant tori on resonant surfaces
 for twist mappings, \href{https://doi.org/10.3934/dcds.2020081}{\textit{Discrete Contin. Dyn. Syst.}} \textbf{40} (2020),
 1389--1409.

\end{thebibliography}
\end{document}